\pdfoutput=1 
\documentclass[reqno]{amsart} 
\usepackage{amssymb,amsmath,amsfonts,amsthm,comment,mathrsfs,times,graphicx}
\usepackage[T1]{fontenc}
\usepackage{lmodern}

\newcommand{\D}{\mathbb{D}}

\newcommand{\Sp}{\mathbb{S}}

\begin{document}

\title[] %
{Not all traces on the circle come from functions of least gradient in the disk}

\author{Gregory S. Spradlin}
\address{Department of Mathematics,Embry-Riddle University, Daytona Beach, Florida 32114-3900, USA}
\email{spradlig@erau.edu}

\author{Alexandru Tamasan}
\address{Department of Mathematics, University of Central Florida, Orlando, Florida, 32816, USA}
\email{tamasan@math.ucf.edu}
\thanks{A.T. supported by NSF, grant DMS-1312883.}


\begin{abstract}
We provide an example of an $L^1$ function on the circle, which cannot be the trace of a function of bounded variation of least gradient in the disk.
\end{abstract}

\maketitle
\numberwithin{equation}{section}
\newtheorem{theorem}{Theorem}[section]
\newtheorem{lemma}[theorem]{Lemma}
\newtheorem{proposition}[theorem]{Proposition}

\section{Introduction}
Sternberg et al. in \cite{SternbergWilliamsZiemer}, and Sternberg and Ziemer in \cite{SternbergZiemer} considered the question of existence, uniqueness and regularity for functions of least gradient and prescribed trace. More precisely for $\Omega\subset\mathbb{R}^n$ a Lipchitz domain, and for a continuous map $g\in C(\partial\Omega)$, they formulate  the problem
\begin{align}\label{min_prob}
\min\{ \int_\Omega|Du|:\;u\in BV(\Omega),\; u|_{\partial\Omega}=g\},
\end{align}
where $BV(\Omega)$ denotes the space of functions of bounded varation, the integral is understood in the sense of the Radon measure $|Du|$ of $\Omega$ and the trace at the boundary is in the sense of the trace of functions of bounded varation.
Solutions to the minimization problem \eqref{min_prob} are called functions of least gradient. For domains $\Omega$ with boundary of non-negative curvature, and which are not locally area minimizing they prove existence, uniqueness and regularity of the solution. Moreover, if the boundary of the domain fails either of the two assumptions they provide counterexamples to existence.

	It is known that traces of functions $f\in BV(\Omega)$ of bounded varation are in $L^1(\partial\Omega)$, and that conversely, any function in $L^1(\partial \Omega)$ admits an extension (in the sense of trace) in $BV(\Omega)$ (in fact in $W^{1,1}(\Omega)$), see e.g., \cite{GiustiBook}. The question we address here is whether solutions of the problem \eqref{min_prob} exist in the case of traces that are merely in 
	$L^1(\partial\Omega)$ and not continuous. We answer this question in the negative by providing a counterexample for the  unit disk, which has a boundary of positive curvature and which is not locally length minimizing.

	Let $\D$  denote the unit disk in the plane and $\Sp$ be its boundary. We prove the following:

\begin{theorem} \label{thm1}
There exists $f \in L^1(\Sp)$ such that the minimization problem
\begin{align}
\min\{ \int_\D |Dw|:\; w \in BV(\D),\; w|_{\Sp} = f\}
\end{align}has no solution.
\end{theorem}

A renewed interest in functions of least gradient with variable weights appeared recently due to its applications to current density impedance imaging, see \cite{NTT11} and references therein. Our counterexample sets a limit on the roughness of the boundary data one can afford to use.

\section{Proof of Theorem~\ref{thm1}}

	We will call the $L^1(\Sp)$-function satisfying 
	Theorem~\ref{thm1} ``$f_\infty$''.  $f_\infty$ is the characteristic function of a fat Cantor set.  Define
		$C_0 \supset C_1 \supset C_2 \supset \cdots$ inductively as follows:
\begin{equation*}
						C_0 = \{(\cos\theta, \sin\theta) \; \Bigm| \;
						        \frac{\pi}{2}-\frac{1}{2} \leq \theta  \leq
										        \frac{\pi}{2}+\frac{1}{2}\},
\end{equation*}
and if $C_n$ consists of $2^n$ disjoint closed arcs, each with arc length
\begin{equation}\label{e2002}
						\theta_n = \frac{1}{2^n} \prod_{i=1}^n (1 - \frac{1}{2^i})
\end{equation}
(if $n = 0$, the ``empty product'' is interpreted as $1$), then $C_{n+1}$ is obtained by removing from the center of each of those arcs an open arc of arc length
$(1 - 1/2^{n+1})\theta_n$. Then $C_{n+1}$ consists of $2^{n+1}$ disjoint closed arcs, each with arc length $\theta_{n+1}$.  For $n = 0, 1, 2, \ldots$, with $\mathcal{H}^1$ denoting one-dimensional Hausdorff measure,
\begin{equation}\label{e2004}
			\mathcal{H}^1(C_n) = 2^n \theta_n = 
			         \prod_{i=1}^{n} (1 - \frac{1}{2^i}) \equiv K_n.
\end{equation}
Define
\begin{equation*}
			C_\infty = \bigcap_{n=0}^\infty C_n.
\end{equation*}
$C_\infty$ is a compact and nowhere dense subset of $\Sp$, with
\begin{equation*}
			\mathcal{H}^1(C_\infty) = \prod_{i=1}^{\infty} (1 - \frac{1}{2^i}) 
			   =\lim_{n \to \infty} K_n \equiv K_\infty > 0.
\end{equation*}
Note that 
$K_\infty$ is well-defined and positive, since all the terms in the infinite product are positive and $\sum_{i=1}^\infty 1/2^i < \infty$.  

	We define $f_\infty\in L^1(\Sp)$ to be the characteristic function of $C_\infty$: 
\begin{equation*}
			f_\infty = \chi_{C_\infty} \in L^1(\Sp).
\end{equation*}
From \cite[Theorem 2.16, Remark 2.17]{GiustiBook} we have that   
\begin{equation*}
  \inf\{\int_\D |Du| \; \Bigm| \;  u \in BV(\D),\ u|_{\Sp} = f_\infty\} \leq
						\|f_\infty\|_{L^1(\Sp)} = K_\infty.
\end{equation*}
We will show that for any $u \in BV(\D)$ with $u|_{\Sp} = f_\infty$,
\begin{equation*}
  \int_\D |Du| > K_\infty,
\end{equation*}
proving Theorem~\ref{thm1}.

\medskip

\begin{figure*}
\centering\includegraphics[scale=1.0]{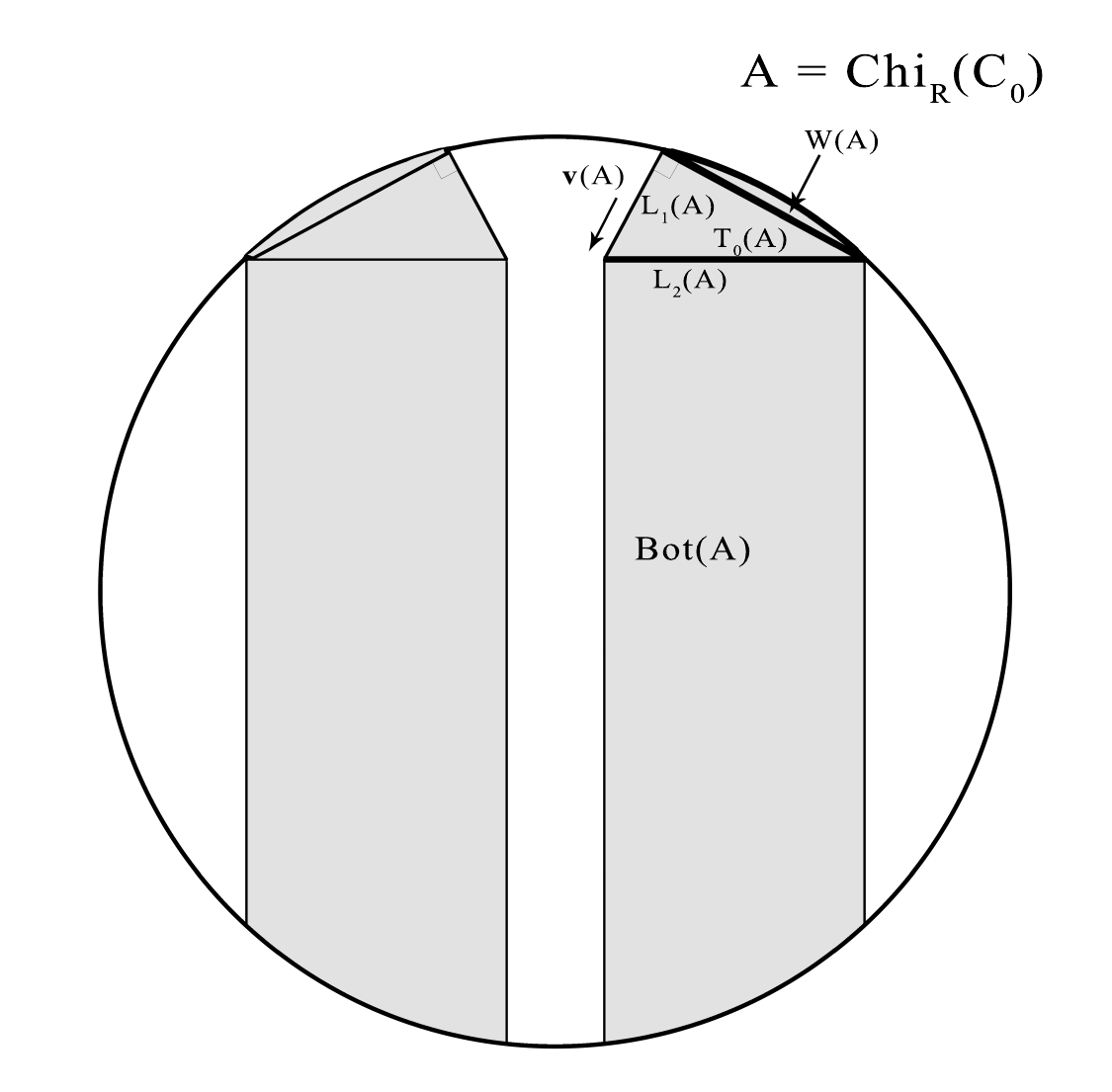}  
\caption{$B_1$} \label{fig:figuretwo}
\end{figure*}

\begin{figure*}
\centering\includegraphics[scale=0.7]{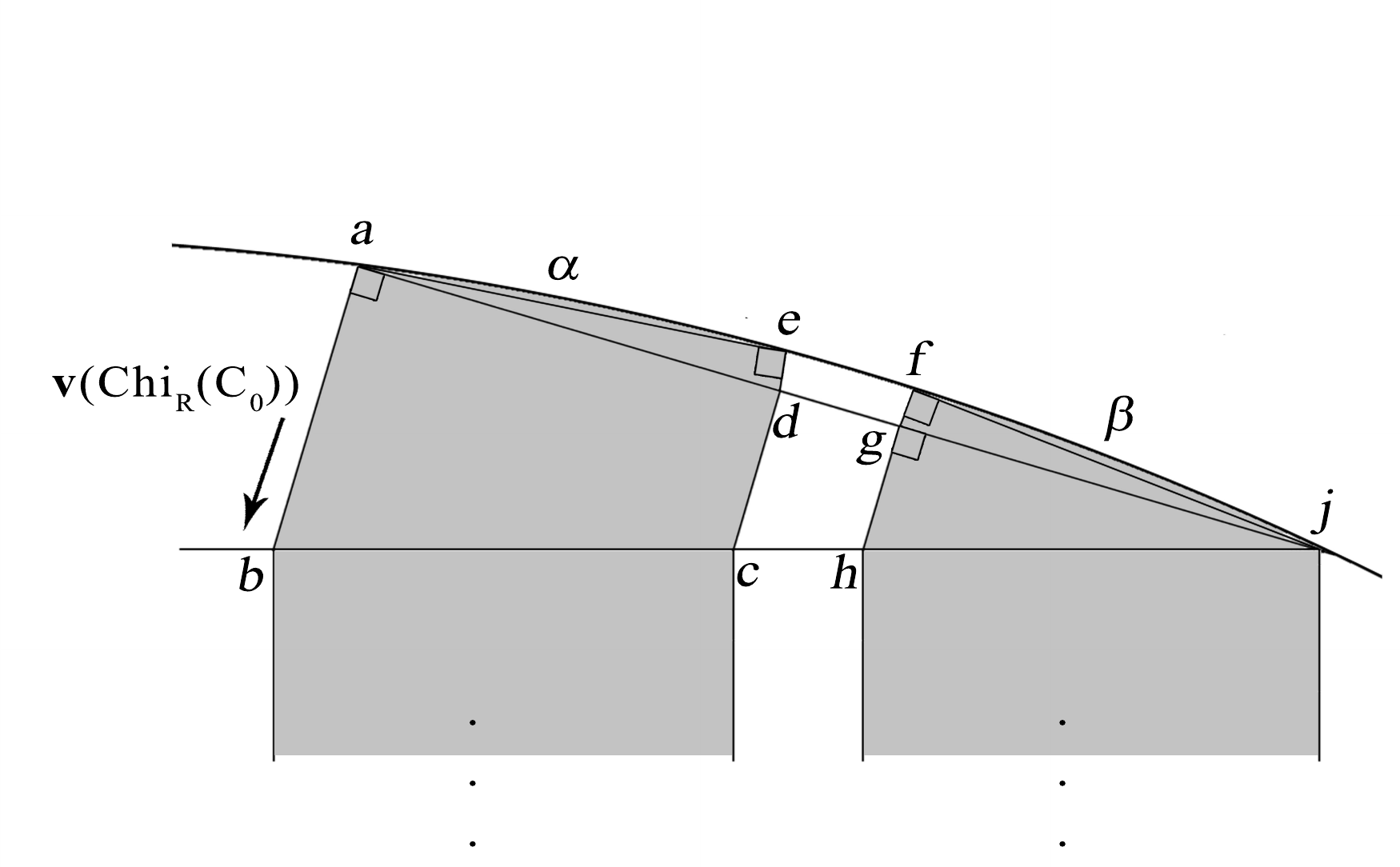} 
\caption{The upper part of the right half of $B_2$.} \label{fig:figurethree}
\end{figure*}

	The idea of the proof is as follows: we construct a compact, nowhere dense subset
$B_\infty$ of ${\overline \D}$ with the property that
\begin{equation}\label{e2010}
\begin{aligned}
&\mbox{(i) If } u \in BV(\D) \mbox{ with } u|_{\Sp} = f_\infty
  \mbox{ and } \int_{\D \setminus B_\infty} |u|\,dx > 0,
                 \mbox{ then } \int_\D |Du| > K_\infty, \\
&\qquad \mbox{ and } \\
&\mbox{(ii) If } u \in BV(\D) \mbox{ with } u|_{\Sp} = f_\infty
   \mbox{ and } \int_{\D \setminus B_\infty} |u|\,dx = 0,
                 \mbox{ then } \int_\D |Du| > K_\infty. \\
\end{aligned}
\end{equation}
Theorem~\ref{thm1} obviously follows from this.  $B_\infty$ has the form
\begin{equation}\label{e2011}
  B_\infty = \bigcap_{n=1}^\infty B_n,
\end{equation}
where $B_1 \supset B_2 \supset B_3 \supset \cdots$, and for each $n \geq 1$, $B_n$ is a compact subset of ${\overline \D}$ with  $2^n$ path components, with each path component the union of a polygon and two circular segments (``circular segment'' is the standard term for the region between an arc and a chord connecting two points on a circle).  That polygon will be defined precisely as the union of at least one triangle with at least one trapezoid. In Figure~\ref{fig:figuretwo}, $B_1$ is the union of the two shaded regions. In Figure~\ref{fig:figurethree}, the two shaded regions constitute the upper portion of the right half of $B_2$.  

	Unfortunately, defining each $B_n$ precisely requires a slew of definitions. For
	$n \geq 0$, $C_n$ is the disjoint union of $2^n$ closed arcs.  Call this collection of arcs $\mathcal{A}_n$.  For example, $\mathcal{A}_0 = \{C_0\}$.  For each
$A \in \mathcal{A}_n$, we will define a set $B_A \subset {\overline \D}$, 
then define $B_n$ as the disjoint union
\begin{equation}\label{e2012}
  B_n = \bigsqcup_{A \in \mathcal{A}_n} B_A.
\end{equation}
Each such $B_A$ is the connected union of a closed circular segment,  $n$ closed polygons, which are all triangles or trapezoids (including at least one triangle), and a ``bottom'' piece that is the union of a trapezoid and a circular segment (in Figure~\ref{fig:figuretwo}, the arc in $\mathcal{A}_1$ in the right half of the $x_1$-$x_2$ plane is called ``$A$'', and $B_A$ is the shaded region in the right half of ${\overline \D}$.  If we call the other arc in $\mathcal{A}_1$ ``$A'$'', then the shaded region in the left half of ${\overline \D}$ is $B_{A'}$.  In Figure~\ref{fig:figurethree}, the shaded region on the left is the top of $B_\alpha$ and the shaded region on the right is the top of $B_\beta$, where $\alpha$ and $\beta$ are the two arcs in $\mathcal{A}_2$ in the right half of the $x_1$-$x_2$ plane.  The other notations used 
in Figures~\ref{fig:figuretwo} 
and \ref{fig:figurethree} will be defined momentarily).
	
	For an arc $A$, let $\operatorname{Cho}(A)$ denote the chord connecting the endpoints of $A$, and $W(A)$ the closed circular segment enclosed by $A$ 
and $\operatorname{Cho}(A)$.  For $A \in \mathcal{A}_n$ with $n \geq 1$, define
$\operatorname{Par}(A)$ (the ``parent'' of $A$) to be the arc in $\mathcal{A}_{n-1}$ containing $A$:
\begin{equation*}
  \operatorname{Par}(A) = A':\ A' \in \mathcal{A}_{n-1},\ A \subset A'.
\end{equation*}
More generally, for $A \in \mathcal{A}_n$ ($n \geq 0$), define
\begin{equation*}
\begin{aligned}
\operatorname{Par}^0(A) &= A,\ \operatorname{Par}^1(A) = \operatorname{Par}(A),
\ \operatorname{Par}^2(A) = \operatorname{Par}(\operatorname{Par}^1(A)), \\
\ \operatorname{Par}^3(A) &= \operatorname{Par}(\operatorname{Par}^2(A)),
\ldots,
\ \operatorname{Par}^n(A) = C_0 \in \mathcal{A}_0.
\end{aligned}
\end{equation*}
For $A \in \mathcal{A}_n$ ($n \geq 0$), define the two ``children'' of $A$,
$\operatorname{Chi}_L(A)$ and $\operatorname{Chi}_R(A)$, by
\begin{equation*}
\begin{aligned}
 &\operatorname{Chi}_L(A), \operatorname{Chi}_R(A) \in \mathcal{A}_{n+1},\
\operatorname{Chi}_L(A), \operatorname{Chi}_R(A) \subset A,\
\operatorname{Chi}_L(A) \cap \operatorname{Chi}_R(A) = \varnothing,\\
&\operatorname{Chi}_L(A) \mbox{ is ``to the left'' or counterclockwise from }
                           \operatorname{Chi}_R(A).\\
\end{aligned}
\end{equation*}
For an arc $A$ of $\Sp$ of arc length less than $\pi$, let $\mathbf{v}(A)$ denote the unit vector perpendicular to $\operatorname{Cho}(A)$ and pointing from $\operatorname{Cho}(A)$ toward $\mathbf{0} \in \mathbb{R}^2$.  For $A \in \mathcal{A}_n$ with $n \geq 1$, let
$T(A)$ denote the unique closed right triangle whose longer leg is
$\operatorname{Cho}(A)$ and whose hypotenuse is a subset of $\operatorname{Cho}(\operatorname{Par}(A))$.

\begin{figure*}
\centering\includegraphics[scale=1.0]{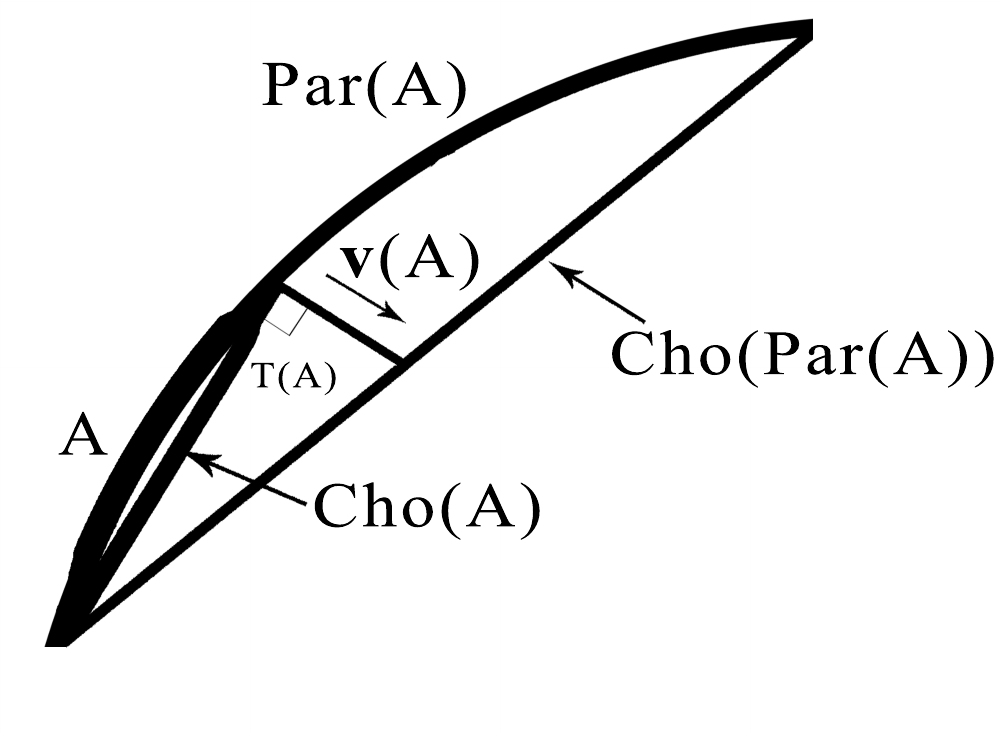}  %
\caption{} \label{fig:figureone}
\end{figure*}

		Figure~\ref{fig:figureone} shows an arc $A$ belonging to $\mathcal{A}_n$ for 
		some $n \geq 1$, along with $\operatorname{Par}(A)$, $\operatorname{Cho}(A)$,
$\operatorname{Cho}(\operatorname{Par}(A))$, $T(A)$, and $\mathbf{v}(A)$.  The lengths of the segments and arcs are not necessarily to scale, and the length of $\mathbf{v}(A)$ is definitely not to scale, since $\mathbf{v}(A)$ is a unit vector and
$\operatorname{Par}(A)$ is an arc of $\Sp$.

		We are finally ready to define $B_A$ 
		(for $A \in \mathcal{A}_n$ with $n \geq 1$).  We will do the $n=1$ and $n=2$ cases first, then the general case.
		
		Suppose $A \in \mathcal{A}_1$ (so $A = \operatorname{Chi}_L(C_0)$ or
		$\operatorname{Chi}_R(C_0)$).  $B_A$ is the union of
$W(A)$, $T(A)$, and a ``bottom'' piece that is the union of a trapezoid and a circular segment.  In order to help establish the pattern for general $n$, we will introduce some notations that are not needed here but will be necessary later.  Define the line segment  $L_1(A) = \operatorname{Cho}(A)$, the triangle $T_0(A) = T(A)$, and define the line segment $L_2(A)$ to be the hypotenuse of $T_0(A)$, which can also be defined
\begin{equation}\label{e2016}
 L_2(A) = \{\mathbf{x} \in \partial T_0(A) 
                 \; \Bigm| \; x_2 = \sin(\frac{\pi}{2} - \frac{1}{2})\}.
\end{equation}
Define the ``bottom'' part of $B_A$, $\operatorname{Bot}(A)$, to be the set of all points in ${\overline \D}$ on or directly ``below'' $L_2(A)$, that is,
\begin{equation}\label{e2017}
 \operatorname{Bot}(A) = \{\mathbf{x} \in {\overline \D} \; \Bigm| \;
          x_2 \leq \sin(\frac{\pi}{2} - \frac{1}{2}),\ x_1 = y_1 \mbox{ for some }
									\mathbf{y} \in L_2(A)\}.
\end{equation}
$\operatorname{Bot}(A)$ is the union of a closed trapezoid and a closed circular segment.  Finally, define
\begin{equation}\label{e2018}
 B_A =W(A) \cup T_0(A) \cup \operatorname{Bot}(A).
\end{equation}
In Lemma~\ref{lemma301} in the Appendix, it is proven that for any $A \in \mathcal{A}_n$
(for $n \geq 0$), $T(\operatorname{Chi}_L(A))$ and $T(\operatorname{Chi}_R(A))$ are disjoint (use $\theta = \theta_n$ and $\alpha = \theta_n/2^{n+1} \geq \theta_n^2/2$, with $\theta_n$
as in \eqref{e2002}).
It follows that $B_{\operatorname{Chi}_L(A)}$ and $B_{\operatorname{Chi}_R(A)}$ are disjoint.
Now $B_1$ is defined as in \eqref{e2012}.  The two shaded regions in 
Figure~\ref{fig:figuretwo} comprise $B_1$.  The arc $\operatorname{Chi}_R(C_0)$ is called $A$, and the parts of $B_A$ (which is the right half of $B_1$) are labelled, along with the vector $\mathbf{v}(A)$, 
which is perpendicular to $\operatorname{Cho}(A)$. The lengths of the segments and arcs are not truly scaled, and
the unit vector $\mathbf{v}(A)$ is drawn with shorter than unit length in order to
fit in the picture.

	Next, suppose $A \in \mathcal{A}_2$.  $B_A$ is the union of a chain of four sets: a closed circular segment, followed by a closed triangle, then a closed triangle or trapezoid, then finally a closed  ``bottom'' piece which is the union of a trapezoid and a circular segment, as in the $n=1$ case.  The intersection of any two consecutive sets in the chain is a line segment.  
	
		Figure~\ref{fig:figurethree} shows the top of the right half of $B_2$, so it shows the top portions of the two rightmost of the four components of $B_2$.  As before,
the lengths of segments and arcs are not necessarily scaled truly, and
$\mathbf{v}(\operatorname{Chi}_R(C_0))$ 
is actually a unit vector, contrary to the picture.  The arc 
$\stackrel{\textstyle\frown}{aj}$ is $\operatorname{Chi}_R(C_0)$. 
For brevity in notation we have defined 
$\alpha = \stackrel{\textstyle\frown}{ae} = 
\operatorname{Chi}_L(\operatorname{Chi}_R(C_0))$ and
$\beta = \stackrel{\textstyle\frown}{fj} = 
\operatorname{Chi}_R(\operatorname{Chi}_R(C_0))$.  The two connected gray regions are the upper portions of $B_\alpha$ and of $B_\beta$.  $B_\alpha$ is the union of 
$W(\alpha)$ (a very thin circular segment in the figure), 
the triangle $\triangle\,ade$, the trapezoid $abcd$, and a ``bottom'' piece 
$\operatorname{Bot}(\alpha)$ consisting of all the
points in ${\overline \D}$ on or directly below the 
line segment ${\overline {bc}}$.  
$B_\beta$ is the union of 
$W(\beta)$ (a very thin circular segment in the figure), the triangle $\triangle\,fgj$, the triangle $\triangle\,ghj$, and a ``bottom'' piece $\operatorname{Bot}(\beta)$ 
consisting of all the points in ${\overline \D}$ 
on or directly below the line segment ${\overline {hj}}$. 

	Generally, for $A \in \mathcal{A}_2$, define the line segment 
$L_1(A) = \operatorname{Cho}(A)$
(so $L_1(\alpha) = {\overline {ae}}$ and $L_1(\beta) = {\overline {fj}}$), and
the triangle $T_0(A) = T(A)$ (so $T_0(\alpha) = \triangle\,ade $ and $T_0(\beta) = \triangle\,fgj$).  Define the line segment
$L_2(A) = {\partial {T_0(A)}} \cap \operatorname{Cho}(\operatorname{Par}(A))$.
$L_2(A)$ can also be described as the hypotenuse of $T_0(A)$.   
In Figure~\ref{fig:figurethree}, 
$L_2(\alpha) = \overline{ad}$ and $L_2(\beta) = \overline{gj}$.
Define $T_1(A) \subset T(\operatorname{Par}^1(A))$ to be the set of all 
points $\mathbf{x}$ in the triangle $T(\operatorname{Par}^1(A))$ with the property that 
for some point $\mathbf{y} \in L_2(A)$,
the vector $\mathbf{x} -\mathbf{y}$ is parallel to 
$\mathbf{v}(\operatorname{Par}^1(A)) \equiv \mathbf{v}(\operatorname{Par}(A))$.  
$T_1(A)$ is either a triangle (this occurs if $A = \operatorname{Chi}_L(\operatorname{Chi}_L(C_0))$ or $\operatorname{Chi}_R(\operatorname{Chi}_R(C_0))$) or a trapezoid 
(this occurs if $A = \operatorname{Chi}_L(\operatorname{Chi}_R(C_0))$ or 
$\operatorname{Chi}_R(\operatorname{Chi}_L(C_0))$).  In Figure~\ref{fig:figurethree}, 
$T_1(\alpha)$ is the trapezoid $abcd$, with
$\alpha = \operatorname{Chi}_L(\operatorname{Chi}_R(C_0))$, and 
$T_1(\beta)$ is the triangle $\triangle\, ghj$, with
$\beta = \operatorname{Chi}_R(\operatorname{Chi}_R(C_0))$.  
$T_1(A)$ can be defined succinctly by
\begin{equation*}
T_1(A) = \{\mathbf{x} \in T(\operatorname{Par}^1(A))   \; \Bigm| \;
       \mathbf{x} -\mathbf{y} \parallel \mathbf{v}(\operatorname{Par}^1(A))
			            \mbox{ for some } \mathbf{y} \in L_2(A)\}.
\end{equation*}
		Define the horizontal line segment $L_3(A)$, similarly to \eqref{e2016}, to be the set of all points in $\partial T_1(A)$ with $x_2$-coordinate $\sin(\pi/2 -1/2)$:		
\begin{equation*}
 L_3(A) = \{\mathbf{x} \in \partial T_1(A) \; \Bigm| \;
                  x_2 = \sin(\frac{\pi}{2} - \frac{1}{2})\}.
\end{equation*}
In other words, $L_3(A)$ is the side of the polygon $\partial T_1(A)$ that is a subset of the horizontal line $\{\mathbf{x} \, \mid \, x_2 = \sin(\pi/2 - 1/2)\}$.  
In Figure~\ref{fig:figurethree}, $L_3(\alpha) = {\overline {bc}}$ 
and $L_3(\beta) = {\overline {hj}}$.  Like in \eqref{e2017}, define the ``bottom'' part of $B_A$, $\operatorname{Bot}(A)$, to be the set of all points in ${\overline \D}$ on or directly below
$L_3(A)$, that is,
\begin{equation*}
 \operatorname{Bot}(A) = \{\mathbf{x} \in {\overline \D} \; \Bigm| \;
                  x_2 \leq \sin(\frac{\pi}{2} - \frac{1}{2}),\, 
									  x_1 = y_1 \mbox{ for some }
									\mathbf{y} \in L_3(A) \}.
\end{equation*}
Like before, $\operatorname{Bot}(A)$ is the union of a trapezoid and a circular segment.  Similar to \eqref{e2018}, define
\begin{equation*}
 B_A = W(A) \cup
     T_0(A) \cup
			 T_1(A) \cup
			 \operatorname{Bot}(A).
\end{equation*}
As in the $n=1$ case, by Lemma~\ref{lemma301} in the Appendix, the sets $B_A$ for the four elements of $\mathcal{A}_2$ are disjoint.  $B_2$ is defined by \eqref{e2012}.  Clearly $B_2 \subset B_1$.  

	Finally we consider the $n > 2$ case.  Let $A \in \mathcal{A}_n$.  $B_A$ is the union of a chain of $n+2$ closed sets: a closed circular segment, followed by $n$ closed polygons which are all triangles or trapezoids, and finally a bottom piece called
$\operatorname{Bot}(A)$ (as before) which is the union of a closed trapezoid and a closed circular segment. Either all $n$ of the polygons are triangles, or 
(more likely), the first $k$ of them are triangles for some $1 \leq k \leq n-1$ 
and
the remaining $n-k$ polygons are trapezoids.  The intersection of any two consecutive sets in the chain is a line segment.
$B_A$ has the form
\begin{equation*}
 B_A = W(A) \cup
			 \bigcup_{k=0}^{n-1} T_k(A) \cup
			 \operatorname{Bot}(A),
\end{equation*}
where $T_k(A)$ and $\operatorname{Bot}(A)$ are defined precisely momentarily.  In order to do so, we must also name the intersections of consecutive sets in the chain, which are line segments, and which we will call $L_1(A), \ldots, L_{n+1}(A)$.
We will also need to use $L_1(A), \ldots, L_{n+1}(A)$ to prove \eqref{e2010}.

		Define
\begin{equation}\label{e2024}
\begin{aligned}
L_1(A) &= \operatorname{Cho}(A),\\
T_0(A) &= T(A),\\
L_2(A) &= \partial T_0(A) \cap \operatorname{Cho}(\operatorname{Par}^1(A)),\\
T_1(A) &= \{\mathbf{x} \in T(\operatorname{Par}^1(A)) \; \Bigm| \;
       \mathbf{x} - \mathbf{y} \parallel \mathbf{v}(\operatorname{Par}^1(A))
			              \mbox{ for some } \mathbf{y} \in L_2(A)\},\\
L_3(A) &= \partial T_1(A) \cap \operatorname{Cho}(\operatorname{Par}^2(A)),\\
T_2(A) &= \{\mathbf{x} \in T(\operatorname{Par}^2(A)) \; \Bigm| \;
       \mathbf{x} - \mathbf{y} \parallel \mathbf{v}(\operatorname{Par}^2(A))
			              \mbox{ for some } \mathbf{y} \in L_3(A)\},\\	
		&\vdots\\
T_{n-1}(A) &= \{\mathbf{x} \in T(\operatorname{Par}^{n-1}(A)) \; \Bigm| \;
       \mathbf{x} - \mathbf{y} \parallel \mathbf{v}(\operatorname{Par}^{n-1}(A))
			              \mbox{ for some } \mathbf{y} \in L_n(A)\},\\
L_{n+1}(A) &= \{\mathbf{x} \in \partial T_{n-1}(A)  \; \Bigm| \;
                   x_2 = \sin(\frac{\pi}{2}-\frac{1}{2})\},\\
\operatorname{Bot}(A) &= \{\mathbf{x} \in \overline{\D}  \; \Bigm| \;
              x_2 \leq \sin(\frac{\pi}{2} - \frac{1}{2}),\ x_1 = y_1 \mbox{ for some }
									\mathbf{y} \in L_{n+1}(A)\}.																		
\end{aligned}
\end{equation}		
Like before, the $B_A$'s are disjoint for all the  $2^n$ arcs $A$ in $\mathcal{A}_n$, and
$B_n$ is defined by \eqref{e2012}.  Clearly $B_1 \supset B_2 \supset B_3 \supset \cdots$.  We define $B_\infty$ by \eqref{e2011}.

	Having defined $B_\infty$, we show that it has property \eqref{e2010}, from which Theorem~\ref{thm1} follows.  This requires three lemmas, followed by an easy proof of \eqref{e2010}(i), then a more involved proof of \eqref{e2010}(ii).
	
\begin{lemma}\label{lemma201} Let $u \in C^\infty(\D) \cap BV(\D)$
with $u|_{\Sp} = f_\infty$, $n \geq 1$, and $A \in \mathcal{A}_n$.  Let $T_0(A), T_1(A), \ldots, T_{n-1}(A)$ and $\operatorname{Bot}(A)$ be as in \eqref{e2024}.  Then
\begin{equation}\label{e2025}
\begin{aligned}
&\int_{W(A)} |\nabla u \cdot \mathbf{v}(A)|\,dx +
     \sum_{k=0}^{n-1}
		    \int_{T_k(A)} |\nabla u \cdot \mathbf{v}(\operatorname{Par}^k(A))|\,dx + \\
				&\qquad +
				   \int_{\operatorname{Bot}(A)} |\nabla u \cdot \mathbf{j}|\,dx \geq
											\cos \left(\frac{K_n}{2^{n+1}}\right)\frac{K_\infty}{2^n}.
\end{aligned}
\end{equation}
\end{lemma}
Here, $\mathbf{j} = \langle 0, 1 \rangle$, as usual.   $K_n$ is from \eqref{e2004}.  There is a slight abuse of notation in \eqref{e2025}: the domain of $u$ is $\D$, not 
${\overline \D}$, but $W(A)$ and $\operatorname{Bot}(A)$ intersect 
$\Sp \equiv \partial \D$, and $T_k(A)$ might intersect $\Sp$.  In all cases, the intersection has $\mathcal{H}^2$-measure zero.  It would be better formally to replace 
``$W(A)$'', ``$T_k(A)$'', and ``$\operatorname{Bot}(A)$'' in \eqref{e2025} with their interiors, or with their intersections with $\D$.  However, this might make the proof of Lemma~\ref{lemma201} less readable, so we will keep the 
notation of \eqref{e2025} in the proof of the lemma, and in the remainder of this section.

	Proof of lemma: define $s_n = 2\sin(K_n/2^{n+1})$, which is the
length of $\operatorname{Cho}(A)$.  Let $L_1(A)$, $\ldots ,$ $ L_{n+1}(A)$ be as
in \eqref{e2024}.  For $k = 1, 2, \ldots, n+1$, let $\phi_k : [0, s_n] \to
{\overline \D}$ be the linear map with $\phi_k(0)$ the 
left endpoint of $L_k(A)$ and $\phi_k(s_n)$ the 
right endpoint of $L_k(A)$ ($L_k(A)$ is not vertical).
Define $\phi_0 :(0,s_n) \to A$ so that $\phi_0(t)$ is the projection of $\phi_1(t)$ onto
$A$ in the direction $-\mathbf{v}(A)$ (the explicit formula for $\phi_0(t)$ is fairly complicated and we do not use it, so we omit it).  
%
%
Now define $g_0, g_1, \ldots, g_{n+1} \in L^1((0,s_n))$ by
\begin{equation*}
		g_0(t) = f_\infty(\phi_0(t)),\ g_k(t) = u(\phi_k(t)) \mbox{ for } 1 \leq k \leq n+1.
\end{equation*}
Now, $\mathcal{H}^1(K_\infty \cap A) = K_\infty/2^n$,
so $\int_A f_\infty \,d\,\mathcal{H}^1 = K_\infty/2^n$.  Recall $\theta_n$ from 
\eqref{e2012}.  Since the angle between $A$ and
$\operatorname{Cho}(A)$ is at most $\theta_n/2 = K_n/2^{n+1}$,
\begin{equation}\label{e2028}
		\|g_0\|_{L^1((0,s_n))} \geq 
		\cos\left(\frac{K_n}{2^{n+1}}\right)\int_A f_\infty \,d\,\mathcal{H}^1 =
		\cos\left(\frac{K_n}{2^{n+1}}\right)\frac{K_\infty}{2^n}.
\end{equation}
Obviously,
\begin{equation*}
	g_0 = (g_0-g_1)+(g_1-g_2)+(g_2-g_3)+\cdots+(g_n-g_{n+1}) + g_{n+1},
\end{equation*}
so by the triangle inequality,
\begin{equation}\label{e2030}
	\|g_0\|_{L^1((0,s_n))} \leq \|g_1-g_0\|_{L^1((0,s_n))} +
	       \sum_{k=1}^n \|g_{k+1}-g_k\|_{L^1((0,s_n))} + \|g_{n+1}\|_{L^1((0,s_n))}.
\end{equation}
Now
\begin{equation}\label{e2031}
	\|g_1-g_0\|_{L^1((0,s_n))} = \int_0^{s_n} |g_1(t) - g_0(t)| \,dt \leq
	    \int_{W(A)} |\nabla u(x) \cdot \mathbf{v}(A)| \,dx.
\end{equation}
For $1\leq k \leq n$, the Fundamental Theorem of Calculus yields
\begin{equation}\label{e2032}
	\|g_{k+1}-g_k\|_{L^1((0,s_n))} = \int_0^{s_n} |g_{k+1}(t) - g_k(t)| \,dt \leq
	   \int_{T_{k-1}(A)}
		    |\nabla u(x) \cdot \mathbf{v}(\operatorname{Par}^{k-1}(A))| \,dx.
\end{equation}
Since $f_\infty = 0$ on the bottom half of $\Sp$, $u|_{\Sp} = f_\infty$, and $\operatorname{Bot}(A) \cap \Sp$ is a subset of the bottom half of
$\Sp$, it follows that
\begin{equation}\label{e2033}
	\|g_{k+1}\|_{L^1((0,s_n))} = \int_0^{s_n} |g_{k+1}(t)| \,dt \leq
						\int_{\operatorname{Bot}(A)} |\nabla u \cdot \mathbf{j}|\,dx.
\end{equation}
Putting \eqref{e2028} and \eqref{e2030}-\eqref{e2033} 
together yields \eqref{e2025}.

\hfill $\square$

	Define $\D_- \subset \D$, the ``lower part'' of $\D$, by
\begin{equation}\label{e2034}
	\D_- = \{\mathbf{x} \in \D \; \Bigm| \; x_2 < \sin(\frac{\pi}{2} - \frac{1}{2})\}.
\end{equation}
From Lemma~\ref{lemma201}, there follows:
\begin{lemma}\label{lemma202} Let $u$ be as in Lemma~\ref{lemma201}: $u \in C^\infty(\D) \cap BV(\D)$ with $u|_{\Sp} = f_\infty$.
Let $n \geq 1$. Then
\begin{equation}\label{e2035}
\begin{aligned}
\sum_{A \in \mathcal{A}_n} &\int_{W(A)} |\nabla u \cdot \mathbf{v}(A)|\,dx +
         \sum_{m=1}^n \sum_{A\in \mathcal{A}_m}
				     \int_{T(A) \cap B_n} |\nabla u \cdot \mathbf{v}(A)|\,dx + \\
		&+ \int_{B_n \cap \D_-} |\nabla u \cdot \mathbf{j}|\,dx \geq
		                    \cos\left(\frac{K_n}{2^{n+1}}\right)K_\infty.
\end{aligned}
\end{equation}
\end{lemma}
Proof: By Lemma~\ref{lemma201},
\begin{equation}\label{e2036}
\begin{aligned}
\sum_{A \in \mathcal{A}_n} &\int_{W(A)} |\nabla u \cdot \mathbf{v}(A)|\,dx +
         \sum_{A \in \mathcal{A}_n} \sum_{k=0}^{n-1}
				     \int_{T_k(A)} |\nabla u \cdot \mathbf{v}(A)|\,dx + \\
		&+ \sum_{A\in \mathcal{A}_n}
		  \int_{\operatorname{Bot}(A)} |\nabla u \cdot \mathbf{j}|\,dx \geq
		                    \cos\left(\frac{K_n}{2^{n+1}}\right)K_\infty.
\end{aligned}
\end{equation}
We must prove that the inequalities \eqref{e2035} and \eqref{e2036} are equivalent.
The right-hand sides and the first terms of the left-hand sides are exactly the same.  
The third summands in the left-hand sides are equal because
$B_n \cap \D_-$ is the disjoint union of the sets $\operatorname{Bot}(A)$ for the $2^n$ arcs $A$ in $\mathcal{A}_n$.  We must show that the second summands on the 
left-hand sides of \eqref{e2035} and \eqref{e2036} are equal.  
Call the common integrand of 
the integrals ``$g(x)$''.  Generally, any two 
distinct sets of the form $T_k(A)$, for $l\geq 1$, 
$A \in \mathcal{A}_l$, and $k\in \{0, \ldots, l-1\}$ have intersection of 
$\mathcal{H}^2$-measure zero.  This includes the case of 
$k=0,\ T_k(A) = T_0(A) \equiv T(A)$.
Therefore the second summands on the left-hand sides 
of \eqref{e2035} and \eqref{e2036} have 
the form $\int_{S_1} g\,dx$ and $\int_{S_2} g\,dx$, where
\begin{equation*}
\begin{aligned}
S_1 &= \bigcup_{m=1}^n \bigcup_{A \in \mathcal{A}_m} \left(T(A) \cap B_n \right) =
    B_n \cap \left( \bigcup_{m=1}^n \bigcup_{A \in \mathcal{A}_m} T(A) \right), \\
S_2 &= \bigcup_{A \in \mathcal{A}_n} \bigcup_{k=0}^{n-1} T_k(A).		
\end{aligned}
\end{equation*}
We will show $S_1 =S_2$ by showing that the two sets are subsets of each other.  First we show  $S_1 \subset S_2$.  Let $m' \in \{1, \ldots, n\}$ and $A' \in \mathcal{A}_{m'}$.  We will show $B_n \cap T(A') \subset S_2$.  Let $A \in \mathcal{A}_n$.  Since 
$T(A') \subset \D \setminus \D_-$ and $\operatorname{Bot}(A) \subset {\overline \D_-}$, 
$x_2 = \sin(\pi/2 -1/2)$ along $T(A') \cap \operatorname{Bot}(A)$.  Now
\begin{equation*}
\operatorname{Bot}(A) \cap 
   \{\mathbf{x} \; \Bigm| \; x_2 = \sin\left(\frac{\pi}{2}-\frac{1}{2}\right)\} =
T_{n-1}(A) \cap 
   \{\mathbf{x} \; \Bigm| \; x_2 = \sin\left(\frac{\pi}{2}-\frac{1}{2}\right)\}.	
\end{equation*}
Thus $T(A') \cap \operatorname{Bot}(A) \subset T_{n-1}(A)$.  
Also, since $m' \leq n$, $T(A') \cap W(A) \subset T(A) \equiv T_0(A)$.  Therefore,
\begin{equation*}
\begin{aligned}
B_n \cap T(A') &\equiv 
    \Big(\bigcup_{A\in \mathcal{A}_n} 
		     \big(W(A) \cup \operatorname{Bot}(A) \cup \bigcup_{k=0}^{n-1} T_k(A)\big) \Big) 
										\cap T(A') = \\
		&= \bigcup_{A\in \mathcal{A}_n} \Big(
		     \big( W(A) \cup \operatorname{Bot}(A) \cup \bigcup_{k=0}^{n-1} T_k(A) \big) 
										\cap T(A') \Big) = \\
	  &= \bigcup_{A\in \mathcal{A}_n} \bigcup_{k=0}^{n-1} (T_k(A) \cap T(A')) \subset
		 \bigcup_{A\in \mathcal{A}_n} \bigcup_{k=0}^{n-1} T_k(A) = S_2,\\
\end{aligned}
\end{equation*}
and $S_1 \subset S_2$.  Next we prove $S_2 \subset S_1$.  Let $A' \in \mathcal{A}_n$ 
and $k' \in \{0, \ldots, n-1\}$.  We will show $T_{k'}(A') \subset S_1$. First,
\begin{equation}\label{e203604}
T_{k'}(A') \subset \bigcup_{k=0}^{n-1} T_k(A') \subset B_{A'} \subset B_n.
\end{equation}
Next, let $m' = n-k' \in \{1, \ldots, n\}$.  
Since $T_{k'}(A') \subset T(\operatorname{Par}^{k'}(A'))$ and 
$\operatorname{Par}^{k'}(A') \in \mathcal{A}_{n-k'} = \mathcal{A}_{m'} $, 
it follows that
\begin{equation}\label{e203605}
T_{k'}(A') \subset T(\operatorname{Par}^{k'}(A')) \subset
             \bigcup_{m=1}^n \bigcup_{A\in \mathcal{A}_m} T(A).
\end{equation}
By \eqref{e203604}, \eqref{e203605}, and the definition of $S_1$,
$T_{k'}(A') \subset S_1$.  
Therefore $S_2 \subset S_1$.  Lemma~\ref{lemma202} is proven.

\hfill $\square$

	Now as $n \to \infty$, $\mathcal{H}^2(\bigcup_{A \in \mathcal{A}_n} W(A)) \to 0$.
Also $K_n \to K_\infty$ as $n \to \infty$. So taking limits of both sides of
\eqref{e2035} as $n \to \infty$ yields the second inequality in the lemma below:
\begin{lemma}\label{lemma203} Let $u \in C^\infty(\D) \cap BV(\D)$ 
with $u|_{\Sp} = f_\infty$.  Then
\begin{equation*}
\int_{B_\infty \cap \D} |\nabla u|\,dx \geq
         \sum_{n=1}^\infty \sum_{A\in \mathcal{A}_n}
				     \int_{T(A) \cap B_\infty} |\nabla u \cdot \mathbf{v}(A)|\,dx
		+ \int_{\D_- \cap B_\infty} |\nabla u \cdot \mathbf{j}|\,dx \geq
		                   K_\infty.
\end{equation*}
\end{lemma}
	The first inequality is obvious because any two different triangles in the 
collection $\{T(A) \mid A \in \mathcal{A}_l,\ l \geq 1\}$ have intersection with 
zero $\mathcal{H}^2$-measure, and all such triangles are disjoint with $\D_-$.

	Now let us prove \eqref{e2010}(i).  Suppose $u \in BV(\D)$ with
$u|_{\Sp} = f_\infty$ and $\int_{\D \setminus B_\infty} |u|\,dx > 0$.  
$\D \setminus B_\infty$ consists of countably many open components.  
For at least one such component $\Omega$, $\int_\Omega |u|\, dx > 0$.  
$\partial \Omega$ contains an arc of $\Sp$ of positive arc length along 
which $f_\infty$ equals zero.  Therefore 
$\int_\Omega |Du| > 0$.  By
\cite[Theorem~1.17, Remark~1.18, Remark~2.12]{GiustiBook}, 
there exists a sequence $(u_m) \subset C^\infty(\D) \cap BV(\D)$ with $u_m|_{\Sp} = f_\infty$ for all $m$, $u_m \to u$ in $L^1(\D)$, 
and $\int_\D |\nabla u_m|\,dx \to \int_\D |Du|$ as $m \to \infty$.  By 
\cite[Theorem~1.19]{GiustiBook}, 
$\lim\inf_{m\to\infty} \int_\Omega |\nabla u_m|\,dx \geq \int_\Omega |Du| > 0$.
Therefore, using Lemma~\ref{lemma203},
\begin{equation*}
\begin{aligned}
\int_\D |Du| &= \lim_{m \to \infty} \int_\D |\nabla u_m|\,dx \geq
         \lim\inf_{m \to \infty} (\int_{\D \cap B_\infty} |\nabla u_m|\,dx  +
					                        \int_\Omega |\nabla u_m|\,dx) \geq \\
	    &\geq B_\infty + \int_\Omega |Du| > B_\infty.
\end{aligned}
\end{equation*}

	Next we prove \eqref{e2010}(ii), which will complete the proof of 
	Theorem~\ref{thm1}.
Suppose $u \in BV(\D)$ with $u|_{\Sp} = f_\infty$ and
$\int_{\D \setminus B_\infty} |u|\,dx = 0$. Recall $\D_-$, defined in \eqref{e2034}. Since $u \neq 0$, $\int_{\D_-} |u|\,dx >0$ or
$\int_{\D \setminus \D_-} |u|\,dx >0$.  We examine the former case first.  Assume
$\int_{\D_-} |u|\,dx >0$.  Then there exists a closed rectangle 
$[a, b] \times [c,d] \subset \D_-$ and $\delta > 0$ with
\begin{equation*}
\int_{[a, b] \times [c,d] } |u|\,dx > \delta.
\end{equation*}
$B_\infty$ is a compact and nowhere dense subset of ${\overline \D}$.  The restriction of $\chi_{B_\infty}$ to $\D_-$ is constant on
vertical line segments.  Therefore there exists an open, dense subset $U$ of $[a,b]$ with
\begin{equation*}
(U \times [c,d]) \cap B_\infty = \varnothing.
\end{equation*}
Let $a < a_1 < b_1 < b$ with $a_1, b_1 \in U$ and
\begin{equation*}
\int_{[a_1, b_1] \times [c,d] } |u|\,dx > \frac{\delta}{2}.
\end{equation*}
Let $(u_m) \subset C^\infty(\D) \cap BV(\D)$ be given by the construction in 
\cite[Theorem~1.17]{GiustiBook}:
$u_m|_{\Sp} = f_\infty$ for all $m$, $u_m \to u$ in $L^1(\D)$ and
$\int_\D |\nabla u|\,dx \to \int_\D |Du| > 0$ as $m \to \infty$.  Furthermore,
the $u_m$'s are obtained by convolving $u$ with $C^\infty$ mollifier functions,
supported on discs, with the radii of the discs approaching $0$ as $m \to \infty$ 
uniformly on the rectangle $[a, b] \times [c,d]$.  Thus, for large enough $m$,
$u_m = 0$ on the vertical line segments $\{a_1\} \times [c,d]$ and
$\{b_1\} \times [c,d]$. By Lemma~\ref{lemma302} in the Appendix,
\begin{equation*}
\int_{[a_1, b_1] \times [c,d] }
    \left|\frac{\partial u_m}{\partial x_1}\right|\,dx \geq
					\frac{2}{b_1 - a_1}\int_{[a_1, b_1] \times [c,d] } |u_m|\,dx >
										\frac{\delta}{b_1-a_1} \equiv \delta_2
\end{equation*}
for large enough $m$.  Clearly, for large enough $m$,
\begin{equation*}
\int_{[a_1, b_1] \times [c,d] }
    \left|\frac{\partial u_m}{\partial x_2}\right|\,dx \leq
				\int_{[a_1, b_1] \times [c,d] } |\nabla u_m|\,dx \leq
				  \int_\D |\nabla u_m|\,dx < 2\int_\Omega |Du|. 								
\end{equation*}
Therefore, for large enough $m$, by Lemma~\ref{lemma303} in the Appendix
(using $g = |\partial u_m/\partial x_1|$ and  
$h = |\partial u_m/\partial x_2| = |\nabla u_m \cdot \mathbf{j}|$),
\begin{equation}\label{e2044}
\begin{aligned}
\int_{[a_1, b_1] \times [c,d] } |\nabla u_m|\,dx &\geq
  \int_{[a_1, b_1] \times [c,d] } |\nabla u_m \cdot \mathbf{j}|\,dx +
	            \frac{\delta_2^2}{4\int_\Omega |Du| + \delta_2} \equiv  \\
				&\equiv \int_{[a_1, b_1] \times [c,d] } |\nabla u_m \cdot \mathbf{j}|\,dx + \delta_3.	
\end{aligned}
\end{equation}
The collection of triangles $\{T(A) \; \Bigm| \; A \in\mathcal{A}_l, \, l \geq 1\}$ is a countable family of sets, for which the intersection of any distinct pair has zero 
$\mathcal{H}^2$-measure. All these triangles are subsets of 
${\overline \D} \setminus \D_-$.  Therefore, applying \eqref{e2044} and
Lemma~\ref{lemma203}, it follows that for large enough $m$,
\begin{equation*}
\begin{aligned}
\int_\D |\nabla u_m| &=
    \int_{\D_-} |\nabla u_m|\,dx  + \int_{\D\setminus \D_-} |\nabla u_m|\,dx   \geq \\
		&\geq  \int_{\D_-} |\nabla u_m \cdot \mathbf{j}|\,dx + \delta_3 +
		        \sum_{n=1}^\infty \sum_{A \in \mathcal{A}_n}
						        \int_{T(A)} |\nabla u_m|\,dx \geq\\
    &\geq \int_{\D_- \cap B_\infty} |\nabla u_m \cdot \mathbf{j}|\,dx + \delta_3 +
		        \sum_{n=1}^\infty \sum_{A \in \mathcal{A}_n}
						        \int_{T(A) \cap B_\infty} |\nabla u_m \cdot \mathbf{v}(A)|\,dx \geq\\
	  &\geq K_\infty + \delta_3.
\end{aligned}
\end{equation*}
Since $\int_\D |\nabla u_m|\,dx \to \int_\D |Du|$ as $m \to \infty$, it follows that
$\int_\D |Du| \geq K_\infty + \delta_3 > K_\infty$.

	Next, suppose $\int_{\D \setminus \D_-} |u|\,dx >0$
(and $\int_{\D \setminus B_\infty} |u|\,dx = 0$). Since
\begin{equation*}
((\D\setminus \D_-)\cap B_\infty) \subset
     \bigcup_{n=1}^\infty \bigcup_{A \in \mathcal{A}_n} T(A),
\end{equation*}
there exists $n' \geq 1$ and $A \in \mathcal{A}_{n'}$ with $\int_{T(A)} |u|\,dx > 0$.
There then exists a closed rectangle $R \subset T(A) \cap \D$ with sides parallel and perpendicular to  $\mathbf{v}(A)$ and $\int_R |u|\,dx > 0$.

	Let $(u_m)$ be given by the construction in 
\cite[Theorem~1.17]{GiustiBook},
as before.  Arguing as before, let the line segment $L$ be one of the 
two sides of $R$
perpendicular to $\mathbf{v}(A)$.  $L$ has an open and 
dense (with respect to the subspace topology on $L$) subset $X$ 
with $X \cap B_\infty = \varnothing$.  From the way
$B_\infty$ is constructed, if $\mathbf{x} \in R$
and the vector $\mathbf{x} - \mathbf{y}$ is parallel to $\mathbf{v}(A)$
for some $\mathbf{y} \in X$, then
$\mathbf{x} \not\in B_\infty$.  Arguing as in the $\int_{\D_-} |u|\,dx > 0$
case, there exists $\delta_3 > 0$ with
\begin{equation}\label{e2047}
  \int_R |\nabla u_m|\,dx \geq \int_R |\nabla u_m \cdot \mathbf{v}(A)|\,dx + \delta_3
\end{equation}
for large enough  $m$.  Using Lemma~\ref{lemma203} 
and \eqref{e2047}, for large enough $m$,
\begin{equation*}
\begin{aligned}
\int_\D |\nabla &u_m|\,dx = \int_{\D_-}|\nabla u_m|\,dx + 
                \int_{\D \setminus \D_-}|\nabla u_m|\,dx \geq \\
			&\geq \int_{\D_- \cap B_\infty} |\nabla u_m|\,dx +
         	 \sum_{n=1}^\infty
						      \sum_{A \in \mathcal{A}_n}
						\int_{T(A)} |\nabla u_m|\,dx \geq \\
			&\geq \int_{\D_- \cap B_\infty} |\nabla u_m \cdot \mathbf{j}|\,dx + \delta_3 +
 \sum_{n=1}^\infty
		\sum_{A \in \mathcal{A}_n}
		\int_{T(A)} |\nabla u_m \cdot \mathbf{v}(A)|\,dx \geq  \\
			&\geq \int_{\D_- \cap B_\infty} |\nabla u_m \cdot \mathbf{j}|\,dx + \delta_3 +
 \sum_{n=1}^\infty
		\sum_{A \in \mathcal{A}_n}
		\int_{T(A) \cap B_\infty} |\nabla u_m \cdot \mathbf{v}(A)|\,dx \geq \\
			&\geq K_\infty +\delta_3.\\
\end{aligned}
\end{equation*}
Like before, since $\int_\D |\nabla u_m|\,dx \to \int_\D |Du|$
as $m \to \infty$, it follows that 
$\int_\D |Du| \geq K_\infty + \delta_3 > K_\infty$.
The proof of Theorem~\ref{thm1} is complete.

\hfill $\square$

\section{Appendix: Three Lemmas}

	This section contains three easy, self-contained lemmas, moved to the end of the paper in order not to interrupt the flow of the main proof.
\begin{lemma}\label{lemma301} Let $\theta \in (0,1]$ and
$\alpha \in [\theta^2/2,\theta)$.  Let $P$ and $S$ be points on $\Sp$ separated by arc length $\theta$.  Let $Q$ and $R$ lie on the arc
$\stackrel{\textstyle\frown}{PS}$,
with
$\stackrel{\textstyle\frown}{QR}$ having arc length $\alpha$,
$\stackrel{\textstyle\frown}{PQ}$ and
$\stackrel{\textstyle\frown}{RS}$
having equal arc length, and $Q$ between $P$ and $R$.  Let $T$ and $U$ lie on the chord
${\overline {PS}}$, chosen such that $\triangle PQT$ and $\triangle RSU$ are right triangles.  Then $\triangle PQT$ and $\triangle RSU$ have disjoint closures.
\end{lemma}
Proof: clearly it suffices to consider $\alpha = \theta^2/2$.  By rotating the arc
 $\stackrel{\textstyle\frown}{PS}$, we may assume $P = (\cos(\theta/2),\sin(\theta/2))$,
$S = (\cos(\theta/2),-\sin(\theta/2))$, $Q = (\cos(\theta^2/4),\sin(\theta^2/4))$, and
$R = (\cos(\theta^2/4),-\sin(\theta^2/4))$.  Define $V = (\cos(\theta/2),0)$.  It suffices to show the angle $\angle PQV$ is obstuse, using a dot product.  We will show
${\vec {QP}} \cdot {\vec {QV}} < 0$.  Using familiar trigonometric identities,
\begin{equation*}
\begin{aligned}
{\vec {QP}} &= \langle \cos(\frac{1}{2}\theta)-\cos(\frac{1}{4}\theta^2),
                       \sin(\frac{1}{2}\theta)-\sin(\frac{1}{4}\theta^2)\rangle,\\
{\vec {QV}} &=\langle \cos(\frac{1}{2}\theta)-\cos(\frac{1}{4}\theta^2),
                     -\sin(\frac{1}{4}\theta^2 \rangle ,\\											
{\vec {QP}} \cdot {\vec {QV}} &= (\cos(\frac{1}{4}\theta^2)-\cos(\frac{1}{2}\theta))^2-
			(\sin(\frac{1}{2}\theta)-\sin(\frac{1}{4}\theta^2))\sin(\frac{1}{4}\theta^2) =\\
			&= \cos^2(\frac{1}{4}\theta^2)+\cos^2(\frac{1}{2}\theta) -
			      2\cos(\frac{1}{4}\theta^2)\cos(\frac{1}{2}\theta) - \\
			&\qquad \qquad \sin(\frac{1}{2}\theta)\sin(\frac{1}{4}\theta^2) +  
		                  \sin^2(\frac{1}{4}\theta^2) =\\				
			&= 1+\frac{1}{2} +\frac{1}{2}\cos(\theta) -
			       (\cos(\frac{1}{2}\theta +\frac{1}{4}\theta^2) + 
						  \cos(\frac{1}{2}\theta -\frac{1}{4}\theta^2)) - \\
			&\qquad \qquad		\frac{1}{2}(\cos(\frac{1}{2}\theta -\frac{1}{4}\theta^2)-
						            \cos(\frac{1}{2}\theta +\frac{1}{4}\theta^2)) =\\
		   &= \frac{3}{2} +\frac{1}{2}\cos\theta-
			       \frac{1}{2}\cos(\frac{1}{2}\theta +\frac{1}{4}\theta^2) -
						   \frac{3}{2}\cos(\frac{1}{2}\theta -\frac{1}{4}\theta^2).
\end{aligned}
\end{equation*}
By the Maclaurin series for $\cos$ and properties of alternating series,
$1 - x^2/2 < \cos x < 1 -x^2/2 + x^4/24$ for $0<x<1$.  
Both $\theta/2+\theta^2/4$ and
$\theta/2-\theta^2/4$ are between $0$ and $1$.  Therefore
\begin{equation*}
\begin{aligned}
{\vec {QP}} \cdot {\vec {QV}} &<
      \frac{3}{2} +\frac{1}{2}(1-\frac{\theta^2}{2}+\frac{\theta^4}{24}) -
			   \frac{1}{2}(1 - \frac{1}{2}(\frac{\theta}{2}+\frac{\theta^2}{4})^2) -
				    \frac{3}{2}(1-\frac{1}{2}(\frac{\theta}{2}-\frac{\theta^2}{4})^2) =\\
				&= -\frac{1}{8}\theta^3 + \frac{1}{24}\theta^4 < 0.
\end{aligned}
\end{equation*}
\hfill 		$\square$

\smallskip

\begin{lemma}\label{lemma302} Let $a<b$, $c<d$, and $u\in C^1([a,b]\times [c,d])$ with
$u=0$ on $\{a,b\} \times [c,d]$.  Then
\begin{equation}\label{e3003}
\int_{y=c}^d \int_{x=a}^b \left|\frac{\partial u}{\partial x}\right|\,dx\,dy \geq
					\frac{2}{b-a}\int_{y=c}^d \int_{x=a}^b |u(x,y)|\,dx\,dy.
\end{equation}
\end{lemma}
Proof: fix $y \in [c,d]$.  Let $x_0 \in (a,b)$ with
$|u(x_0,y)| = \max_{[a,b]\times \{y\}} |u|$.  Then
\begin{equation}\label{e3004}
\begin{aligned}
\int_{x=a}^b |u(x,&y)|\,dx \leq (b-a)|u(x_0,y)| = \\
   &= (\frac{b-a}{2})(\left|u(x_0,y)-u(a,y)\right| + \left|u(b,y)-u(x_0,y)\right|) = \\
	&= (\frac{b-a}{2})(\left|\int_a^{x_0} \frac{\partial u}{\partial x}\,dx\right| +
	       \left|\int_{x^0}^b \frac{\partial u}{\partial x}\,dx\right|) \leq \\
	  &\leq (\frac{b-a}{2})(\int_a^{x_0} \left|\frac{\partial u}{\partial x}\right|\,dx +
		            \int_{x_0}^b \left|\frac{\partial u}{\partial x}\right|\,dx) =
		(\frac{b-a}{2})\int_a^b \left|\frac{\partial u}{\partial x}\right|\,dx.
\end{aligned}
\end{equation}
Multiplying both sides of \eqref{e3004} by $2/(b-a)$ and integrating from $y=c$ to 
$y=d$ yields \eqref{e3003}.

\hfill $\square$

\smallskip

\begin{lemma}\label{lemma303} Let $M, \delta > 0$, let $\Omega$ be an open subset of
$\mathbb{R}^n$ ($n \geq 1$), let $g, h \in L^1(\Omega)$ 
with $g, h \geq 0$ Lebesgue-a.e., and assume
$\int_\Omega g\,dx \geq \delta$, $\int_\Omega h\,dx \leq M$.  Then
\begin{equation}\label{e3005}
\int_\Omega \sqrt{g^2 + h^2} \,dx \geq \int_\Omega h \, dx +
                        \frac{\delta^2}{2M+\delta}.
\end{equation}
\end{lemma}
Proof (this proof is courtesy of Oleksiy Klurman of the University of Manitoba): 
since $x^2/(\sqrt{x^2+y^2} + |y|) \to 0$ as $(x,y) \to (0,0)$,
we will interpret the expression ``$g^2/(\sqrt{g^2+h^2} + h)$'' 
as zero when $g=h=0$ below.
By the Cauchy-Schwarz Inequality,
\begin{equation*}
\begin{aligned}
(\int_\Omega g\,dx)^2 &=
    \left(\int_\Omega \frac{g}{\sqrt{\sqrt{g^2+h^2}+h}} \cdot
		               {\sqrt{\sqrt{g^2+h^2}+h}}\, dx \right)^2 \leq \\
	&\leq \left(\int_\Omega \frac{g^2}{\sqrt{g^2+h^2}+h}\,dx \right)
	      \left(\int_\Omega {\sqrt{g^2+h^2}+h}\,dx \right) =\\
	&= \left(\int_\Omega {\sqrt{g^2+h^2} -h}\,dx \right)
	   \left(\int_\Omega {\sqrt{g^2+h^2}+h}\,dx \right) \leq \\
	&\leq \left(\int_\Omega {\sqrt{g^2+h^2} -h}\,dx \right)
	      \left(2\int_\Omega h\,dx + \int_\Omega g\,dx \right) \leq \\
	&\leq \left(\int_\Omega {\sqrt{g^2+h^2} -h}\,dx \right)
	      \left(2M + \int_\Omega g\,dx \right).
\end{aligned}
\end{equation*}
Therefore,
\begin{equation}\label{e3007}
\int_\Omega \sqrt{g^2 +h^2} -h\,dx \geq
         \frac{(\int_\Omega g\,dx)^2}{2M + \int_\Omega g \,dx} \geq
				              \frac{\delta^2}{2M + \delta},
\end{equation}
because the map $x \mapsto x^2/(2M + x)$ is an increasing function of $x$ for
$x \geq 0$.  Rearranging \eqref{e3007} yields \eqref{e3005}.

\hfill $\square$


\vfill

\pagebreak
{\bf Acknowledgment:} {The authors would like to thank Christina Spradlin for producing the images used in this document.}

\end{document}